\documentclass[11pt]{article}
\usepackage{amssymb}
\usepackage{times}

\title{Classes de Wadge potentielles des bor\'eliens \`a coupes d\'enombrables.\indent}
\author{Dominique LECOMTE}
\date{\it ~C. R. Acad. Sci. Paris~\rm 317, S\'erie 1 (1993), 1045-1048}

\newcommand{\Ana}{{\it\Sigma}^{1}_{1}}

\newcommand{\boraone}{{\bf\Sigma}^{0}_{1}}
\newcommand{\boratwo}{{\bf\Sigma}^{0}_{2}}
\newcommand{\borathree}{{\bf\Sigma}^{0}_{3}}

\newcommand{\boraxi}{{\bf\Sigma}^{0}_{\xi}}

\newcommand{\bortwo}{{\bf\Delta}^{0}_{2}}
\newcommand{\borone}{{\bf\Delta}^{0}_{1}}

\newcommand{\bormone}{{\bf\Pi}^{0}_{1}}
\newcommand{\bormtwo}{{\bf\Pi}^{0}_{2}}
\newcommand{\bormthree}{{\bf\Pi}^{0}_{3}}

\newtheorem{thm} {Th\'eor\`eme} [section]

\newtheorem{defis} [thm] {D\'efinitions}
\newtheorem{cor} [thm] {Corollaire}

\newtheorem{prop} [thm] {Proposition}

\begin{document}

\maketitle

\noindent {\footnotesize {\bf R\'esum\'e.} On donne, pour chaque classe de Wadge non auto-duale 
$\Gamma$ contenue dans la classe des $G_\delta$, une caract\'erisation des 
bor\'eliens qui ne sont pas potentiellement dans $\Gamma$, parmi les bor\'eliens \`a  
coupes verticales d\'enombrables ; pour ce faire, on utilise des r\'esultats 
d'uniformisation partielle.}\bigskip

\centerline{\bf Potential Wadge classes of Borel sets with countable sections.}\bigskip

\noindent {\footnotesize {\bf Abstract.} We give, for each non self-dual Wadge class $\Gamma$ contained 
in the class of the $G_\delta$ sets, a characterization of Borel sets which are not 
potentially in $\Gamma$, among Borel sets with countable vertical sections; 
to do this, we use results of partial uniformization.}\bigskip\smallskip

 On va poursuivre dans cette Note l'\'etude des classes de Wadge potentielles 
entam\'ee dans [L]. On utilisera les notations standard de la th\'eorie
descriptive des ensembles, qui peuvent \^etre trouv\'ees
dans [Mo]. Par exemple, si $\xi$ est un ordinal d\'enombrable impair, 
on notera $D_\xi(\boraone)$ la classe des ensembles de la forme 
$\bigcup_{\eta<\xi,~\eta\mbox{~pair}}$ $U_\eta\setminus 
(\bigcup_{\theta<\eta} U_\theta)$, o\`u $(U_\eta)_{\eta<\xi}$ est une suite 
croissante d'ouverts. Rappelons les d\'efinitions de base :\bigskip

\noindent\bf D\'efinitions.\it\ (a) Soit $\Gamma$ une classe de parties d'espaces 
polonais de dimension 0. On dit que $\Gamma$ est une $classe\ de\ Wadge$ s'il existe 
un espace polonais $P_0$ de dimension 0, et un bor\'elien $A_0$ de 
$P_0$ tels que pour tout espace polonais $P$ de dimension 0 et pour toute partie 
$A$ de $P$, $A$ est dans $\Gamma$ si et seulement s'il existe une 
fonction continue $f$ de $P$ dans $P_0$ telle que $A = f^{-1}(A_0)$.\smallskip

\noindent (b) Soient $X$ et $Y$ des espaces polonais, et $A$ un bor\'elien de 
$X \times Y$. Si $\Gamma$ est une classe de Wadge, on dira que 
$A$ est $potentiellement\ dans\ \Gamma$ $($ce qu'on notera $A \in\mbox{pot}(\Gamma))$ 
s'il existe des topologies polonaises de dimension 0, $\sigma$ $($sur $X)$ et 
$\tau$ $($sur $Y)$, plus fines que les topologies initiales, telles que $A$, 
consid\'er\'e comme partie de $(X, \sigma) \times (Y, \tau)$, soit dans $\Gamma$.\rm\bigskip

 Dans l'\'etude des relations d'\'equivalence bor\'eliennes, par exemple dans 
[HKL], on \'etudie le pr\'e-ordre qui suit. Si $E$ $($resp. $E')$ est une 
relation d'\'equivalence bor\'elienne sur l'espace polonais $X$ $($resp. $X')$, on 
pose
$$E \leq E' \Leftrightarrow\ \mbox{il existe }f\mbox{ bor\'elienne de }X\mbox{ dans }X'\mbox{ telle 
que }x E y\Leftrightarrow f(x)E' f(y).$$
La derni\`ere relation peut s'\'ecrire $E = (f\times f)^{-1} (E')$ ; or si 
$E'$ est dans $\Gamma$ $($ou m\^eme si $E'$ est $\mbox{pot}(\Gamma))$ et $E = (f\times 
f)^{-1} (E')$, $E$ est $\mbox{pot}(\Gamma)$. Ceci motive l'introduction de la 
notion de classe de Wadge potentielle.

\vfill\eject

 Les classes de Wadge envisag\'ees dans 
cette Note seront les classes de Baire $\boraxi$, les $D_\xi(\boraone)$, et 
leurs classes duales (la classe duale $\check \Gamma$ de la classe $\Gamma$ 
est la classe des compl\'ementaires des \'el\'ements de $\Gamma$ ; par exemple, 
$\check \boraone = \bormone$). L'article [L] sugg\'erait que l'\'etude des probl\`emes 
d'uniformisation partielle pourrait \^etre riche d'enseignements pour l'\'etude 
des classes de Wadge potentielles ; cet espoir est confirm\'e par cette Note.

\section{$\!\!\!\!\!\!$ Quelques r\'esultats sur l'uniformisation partielle.}\indent

 On commence par donner un nouvel exemple de th\'eor\`eme vrai pour la mesure 
et pas pour la cat\'egorie. Mauldin a d\'emontr\'e dans [Ma] le r\'esultat suivant :

\begin{thm} Soient $X$ et $Y$ des espaces 
polonais, $\lambda$ (resp. $\mu$) une mesure de 
probabilit\'e sur $X$ (resp. $Y$), et $A$ un bor\'elien de 
$X\times Y$ ayant ses coupes horizontales (resp. 
verticales) non d\'enombrables $\mu$-presque partout (resp. 
$\lambda$-presque partout). Alors il existe un bor\'elien 
$F$ de $X$ et un bor\'elien $G$ de $Y$ tels que $\lambda (F) = 
\mu (G) = 1$, et un isomorphisme bor\'elien de $F$ sur $G$ dont 
le graphe est contenu dans $A$.\end{thm}

 On peut se demander si on a un r\'esultat analogue en 
rempla\c cant ``ensemble de mesure nulle" par ``ensemble maigre". 
On va voir que non.

\begin{defis} (a) Un $G_{\delta}$ d'un espace topologique est 
dit $presque\mbox{-}ouvert$ (ou p.o.) s'il est contenu dans l'int\'erieur de son 
adh\'erence (ce qui revient \`a dire qu'il est dense dans un ouvert).\smallskip

\noindent (b) Si $X$ et $Y$ sont des espaces topologiques, une partie $A$
de $X\times Y$ sera dite $localement\ \grave a\ projec\mbox{-}$ $tions\ ouvertes$ (ou l.p.o.) si 
pour tout ouvert $U$ de $X \times Y$, les projections de $A\cap U$ sont 
ouvertes.\end{defis}

 Les ensembles l.p.o. se rencontrent par exemple dans la situation suivante : 
$A$ est $\Ana$ dans un produit de deux espaces polonais r\'ecursivement 
pr\'esent\'es. Si on munit ces deux espaces de leur topologie de 
Gandy-Harrington (celle engendr\'ee par les $\Ana$), $A$ devient l.p.o. dans le 
nouveau produit. C'est essentiellement dans cette situation qu'on utilisera 
cette notion, au cours des preuves dans la section 2. Posons 
$$A_0 := \{ (x,K) \in 2^\omega\times {\cal K}(2^\omega)\setminus 
\{ \emptyset \} ~/~x\in K \}\mbox{,}$$ 
et soient $F\subseteq2^\omega$, et $G\subseteq {\cal K}(2^\omega)\setminus \{ \emptyset \}$.

\begin{thm} Le ferm\'e $A_0$ est l.p.o. \`a coupes 
verticales non d\'enombrables, a ses coupes horizontales non d\'enombrables 
sur un ensemble co-maigre, et \smallskip

\noindent  (a) Si $F$ et $G$ sont co-maigres et $f:F\rightarrow G$ 
est un isomorphisme bor\'elien, $\mbox{Gr}(f) \not \subseteq A_0$.\smallskip

\noindent (b) Si $F$ et $G$ sont presque-ouverts non vides et 
$f:F\rightarrow G$ surjective continue ouverte, $\mbox{Gr}(f) \not 
\subseteq A_0$.\end{thm}

 Il r\'esulte imm\'ediatement du th\'eor\`eme 19.6 de [O] 
que sous l'hypoth\`ese du continu, il existe une bijection 
idempotente $\Phi$ de $[0,1]$ sur lui-m\^eme telle que $E$ 
est maigre ssi ${\lambda (\Phi^{-1}(E)) = 0}$, o\`u 
$\lambda$ est la mesure de Lebesgue sur $[0,1]$. En analysant les 
raisons du r\'esultat n\'egatif 1.3, on peut montrer le

\begin{cor} Une telle fonction $\Phi$ n'est 
pas bor\'elienne. De plus, sous l'hypoth\`ese de d\'etermination 
des jeux ${\bf\Delta}^1_{2n+3}$,~\rm $\Phi$ n'est pas 
${\bf\Pi}^1_{2n+1}\rm$-mesurable.\end{cor}

\vfill\eject
 
 Le r\'esultat essentiel de cette section est le 

\begin{thm} Soient $X$ et $Y$ des espaces 
polonais parfaits de dimension 0, $A$ un $G_{\delta}$ l.p.o. 
non vide de $X\times Y$. Alors il existe des p.o. non vides $F$ et $G$ 
(l'un contenu dans $X$, et l'autre dans $Y$), et $f:F\rightarrow G$ 
surjective continue ouverte dont le graphe est contenu dans $A$, ou dans l'ensemble 
$A^* := \{ (y,x)~/~(x,y) \in A \}$.\end{thm}

 La partie (b) du th\'eor\`eme 1.3 montre qu'on ne peut pas, en 
g\'en\'eral, avoir uniformisation dans les deux sens, malgr\'e des hypoth\`eses 
sym\'etriques. Cependant, on a uniformisation si l'on n'exige pas que l'image soit ``grosse" :

\begin{prop} Soient $X$ et $Y$ des espaces 
m\'etrisables s\'eparables de dimension 0, $Y$ \'etant complet, et 
$f:Y\rightarrow X$ une surjection continue ouverte ; alors 
il existe un hom\'eomorphisme $g$ de $X$ sur une partie de $Y$ tel que
 $f \circ g =\mbox{Id}_X$.\end{prop}

 Bien s\^ur, ce r\'esultat est faux si on enl\`eve la condition de dimension sur $X$.

\begin{cor} Soient $X$ et $Y$ des espaces 
polonais parfaits de dimension 0, $A$ un $G_{\delta}$ l.p.o. 
non vide de $X\times Y$. Alors $A$ est uniformisable sur un 
presque-ouvert non vide de $X$ par une application continue 
et ouverte sur son image.\end{cor}
 
 Cette image est en g\'en\'eral rare, d'apr\`es ce qui pr\'ec\`ede. 

\section{$\!\!\!\!\!\!$ Applications aux classes de Wadge potentielles.}\indent

 Soit $\xi$ un ordinal d\'enombrable non nul. On d\'efinit 
$f:\omega^{<\omega} \rightarrow \{-1\} \cup (\xi +1)$, par 
r\'ecurrence sur $\vert s\vert $, comme suit : 
$f(\emptyset) = \xi$ et
$$f(s^\frown n) = \left\{\!\!\!\!\!\!
\begin{array}{ll} 
& \bullet~-1~\mbox{si}~f(s)\leq 0\mbox{,}\cr 
& \bullet~\theta~\mbox{si}~f(s)=\theta +1\mbox{,}\cr
& \bullet~\mbox{un~ordinal~impair~de}~f(s)~\mbox{tel~que~la~suite}~
(f(s^\frown n))_n~\mbox{soit~co-finale~dans} \cr & f(s)~\mbox{et~
strictement~croissante~si}~f(s)~\mbox{est~limite~non~nul.}
\end{array}
\right.$$
On d\'efinit alors des arbres : 
$T_{\xi} := \{ sÊ\in \omega^{<\omega}~/~f(s) \not= -1 \}$ et 
$T'_{\xi} := \{ sÊ\in T_{\xi}~/~f(s) \not= 0 \}$. Dans la suite, si $f_s$ est une fonction partielle de $X$ 
dans $Y$ ou de $Y$ dans $X$, on notera $G_s$ la partie de 
$X\times Y$ \'egale au graphe de $f_s$ si $f_s$ va de $X$ 
dans $Y$, et \'egale \`a $\mbox{Gr}(f_s)^*$ sinon.

\begin{thm} Soient $X$ et $Y$ des espaces 
polonais, $A$ un bor\'elien $\mbox{pot}(\borathree)$ et $\mbox{pot}(\bormthree)$ de 
$X \times Y$, et $\xi$ un ordinal d\'enombrable non nul.\smallskip

\noindent (a) Si $\xi$ est pair, $A$ est non-$\mbox{pot}(D_{\xi}(\boraone))$ 
si et seulement s'il existe des espaces polonais parfaits 
$Z$ et $T$ de dimension 0, des ouverts-ferm\'es non vides 
$A_s$ et $B_s$ (l'un dans $Z$ et l'autre dans $T$, pour 
$s$ dans $T_{\xi}$), des surjections continues ouvertes 
$f_s$ de $A_s$ sur $B_s$, et des injections continues $u$ 
et $v$ tels que si 
$B_p := \bigcup_{s \in T_{\xi}~/~
\vert s\vert ~\mbox{paire}} G_s$ et $B_i := \bigcup_{s \in T_{\xi}~/~
\vert s\vert ~\mbox{impaire}} G_s$, on ait $\overline{B_p} = B_p \cup B_i$, 
$B_p \subseteq (u\times v)^{-1}(A)$, $B_i \subseteq (u\times v)^{-1}(\check A)$, 
et $G_s = \overline{\bigcup_{n \in \omega} G_{s^\frown n}} \setminus 
(\bigcup_{n \in \omega} G_{s^\frown n})$ si $s \in T'_{\xi}$.\smallskip

\noindent (b) Si $\xi$ est impair, $A$ est non-$\mbox{pot}(\check D_{\xi}(\boraone))$ si et 
seulement s'il existe des espaces polonais $Z$ 
et $T$ parfaits de dimension 0, des ouverts-ferm\'es non vides $A_s$ 
et $B_s$ (l'un dans $Z$ et l'autre dans $T$, pour $s$ dans 
$T_{\xi}$), des surjections continues ouvertes 
$f_s$ de $A_s$ sur $B_s$, et des injections continues $u$ 
et $v$ tels que si $B_p := \bigcup_{s \in T_{\xi}~/~
\vert s\vert ~\mbox{paire}} G_s$ et $B_i := \bigcup_{s \in T_{\xi}~/~
\vert s\vert ~\mbox{impaire}} G_s$, on ait $\overline{B_i} = B_p \cup B_i$, 
$B_i \subseteq (u\times v)^{-1}(A)$, $B_p \subseteq (u\times v)^{-1}(\check A)$, 
et $G_s = \overline{\bigcup_{n \in \omega} G_{s^\frown n}} \setminus 
(\bigcup_{n \in \omega} G_{s^\frown n})$ si $s \in T'_{\xi}$.\end{thm} 
 
 L'hypoth\`ese ``$A\in \mbox{pot}(\borathree)\cap \mbox{pot}(\bormthree)$" recouvre en particulier le cas o\`u $A$ est \`a coupes verticales d\'enombrables, ou \`a coupes verticales co-d\'enombrables.
 
\begin{thm} Soient $X$ et $Y$ des espaces 
polonais, $A$ un bor\'elien \`a coupes verticales d\'enombrables de 
$X \times Y$. \smallskip

\noindent (a) $A$ est non-$\mbox{pot}(\bormone)$ si et seulement s'il 
existe des espaces polonais $Z'$ et $T'$ parfaits de dimension 0, 
une suite d'ouverts-ferm\'es $(A_n)$ $($resp. $(B_n))$ de $Z'$ $($resp. $T')$, 
des surjections continues ouvertes $f_n$ de $A_n$ sur $B_n$, et des 
fonctions continues $U$ et $V$ tels que si $B = \bigcup_{n>0} \mbox{Gr}(f_n)$, 
on ait les inclusions $\emptyset \not= \mbox{Gr}(f_0) \subseteq \overline{B} \setminus B$ et 
$B=(U\times V)^{-1}(A)$.\smallskip

\noindent (b) $A$ est non-$\mbox{pot}(\bormone)$ si et seulement s'il 
existe des espaces polonais $Z$ et $T$ parfaits de dimension 0 non vides, 
une suite d'ouverts denses $(E_n)$ de $Z$, 
des applications continues ouvertes $g_n$ de $E_n$ dans $T$, et des 
fonctions continues $u$ et $v$ tels que $(g_n)_{n>0}$ converge simplement 
vers $g_0$ sur $\bigcap_{n\in\omega} E_n$, $\mbox{Gr}(g_0) \subseteq 
(u\times v)^{-1}(\check A)$ et 
$\bigcup_{n>0} \mbox{Gr}(g_n) \subseteq (u\times v)^{-1}(A)$.\end{thm}

 On note, si $\Gamma$ est une classe de Wadge non auto-duale, $\Gamma^+$ la 
classe de Wadge successeur de $\Gamma$ pour l'inclusion.

\begin{prop} Soient $X$ et $Y$ des espaces polonais, $A$ un 
bor\'elien \`a coupes verticales d\'enombra-bles de $X \times Y$, et $k$ un entier 
naturel non nul. Alors on a les \'equivalences suivantes :\smallskip

\noindent (a) $A$ est non-$\mbox{pot}(\borone)\Leftrightarrow A$ est non-$\mbox{pot}(\boraone)
\Leftrightarrow$ la projection de $A$ sur $Y$ est non 
d\'enombrable.\smallskip

\noindent (b) $A$ est non-$\mbox{pot}(\check D_{2k-1}(\boraone))\Leftrightarrow A$ 
est non-$\mbox{pot}(\check D_{2k}(\boraone))\Leftrightarrow A$ est 
non-$\mbox{pot}(D_{2k-1}(\boraone)^+)$.\smallskip

\noindent (c) $A$ est non-$\mbox{pot}(D_{2k+1}(\boraone))\Leftrightarrow A$ est 
non-$\mbox{pot}(D_{2k}(\boraone))\Leftrightarrow A$ est 
non-$\mbox{pot}(D_{2k}(\boraone)^+)$.\end{prop} 
  
 Les bor\'eliens \`a coupes verticales d\'enombrables sont $\mbox{pot}(\boratwo)$, donc il 
nous reste \`a caract\'eriser lesquels de ces bor\'eliens sont $\mbox{pot}(\bormtwo)$ (en 
effet, les seules classes de Wadge non auto-duales contenues dans 
$\bortwo$ sont $\{\emptyset\}$, $D_\xi(\boraone)$, et leurs classes duales). 
On a le r\'esultat classique d'Hurewicz, d\'emontr\'e dans [SR] :

\begin{thm} Soit $X$ un espace polonais, et $A$ un bor\'elien de $X$. 
Alors $A$ est non-$\bormtwo$ si et seulement s'il existe $E$ d\'enombrable 
sans point isol\'e tel que $\overline{E}\setminus E \approx \omega^\omega$ et $E = A \cap 
\overline{E}$.\end{thm}
 
\begin{thm} Soient $X$ et $Y$ des espaces polonais, 
$A$ un bor\'elien de $X\times Y$ \`a coupes verticales 
d\'enombrables. Alors $A$ est non-$\mbox{pot}(\bormtwo)$ si et 
seulement s'il existe des espaces polonais $Z$ et $T$ parfaits de dimension
 0 non-vides, des injections continues $u$ et 
$v$, des ouverts denses $(A_n)$ de $Z$, des applications 
continues et ouvertes $f_n$ de $A_n$ dans $T$, tels que pour tout $x$ 
dans $\bigcap_{n\in\omega} A_n$, l'ensemble $E_x := \{ 
f_n(x)~/~n\in\omega \}$ soit sans point isol\'e, 
$\overline{E_x} \setminus E_x \approx \omega^\omega$, et 
$E_x = (u\times v)^{-1}(A)_x \cap \overline{E_x}$.\end{thm}

\section{$\!\!\!\!\!\!$ R\'ef\'erences.}

\noindent [HKL]\ \ L. A. Harrington, A. S. Kechris et A. Louveau,~\it A Glimm-Effros 
dichotomy for Borel equivalence relations,\rm~J. Amer. Math. Soc.~3 (1990), 903-928

\noindent [L]\ \ D. Lecomte,~\it Classes de Wadge potentielles et th\'eor\`emes d'uniformisation 
partielle,~\rm Fund. Math. 143 (1993), 231-258

\noindent [Ma]\ \ R. D. Mauldin,~\it One-to-one selections, marriage theorems,~\rm Amer. J. Math.~104 (1982), 823-828

\noindent [Mo]\ \ Y. N. Moschovakis,~\it Descriptive set theory,~\rm North-Holland, 1980

\noindent [O]\ \ J. C. Oxtoby,~\it Measure and category,~\rm Springer-Verlag, 1971

\noindent [SR]\ \ J. Saint Raymond,~\it La structure bor\'elienne d'Effros est-elle standard ?,~\rm 
Fund. Math.~100 (1978), 201-210

\end{document}